\documentclass[graybox]{svmult}

\usepackage{type1cm}        
\usepackage{makeidx}         
\usepackage{graphicx}        
\usepackage{multicol}        
\usepackage[bottom]{footmisc}

\usepackage{newtxtext}       %
\usepackage{newtxmath}       

\usepackage{stmaryrd} 

\makeindex             

\begin{document}

\newcommand{\eps}{\varepsilon} \renewcommand{\det}{\mathrm{det}} \newcommand{\argmin}{ \mathrm{argmin} \,}
\newcommand{\Om}{\Omega} \def\interior{\mathaccent'27} 
\newcommand{\weakto}{ \rightharpoonup}  \newcommand{\weakstarto}{\stackrel{*}{\rightharpoonup}}
\newcommand{\R}{\mathbb{R}}
\newcommand{\stress}{\boldsymbol{\sigma}} \newcommand{\strain}{\boldsymbol{\epsilon}} 
\newcommand{\neu}{ \partial_{\mbox{\it \tiny N}} } \newcommand{\dir}{\partial_{\mbox{\it \tiny D}} }
\newcommand{\jump}[1]{\llbracket #1 \rrbracket}
\newcommand{\separe}{\medskip}
\newcommand{\F}{\mathcal{F}}\renewcommand{\E}{\mathcal{E}}


\title*{A quasi-static model for craquelure patterns}
\author{M. Negri}
\institute{Department of Mathematics - University of Pavia \email{matteo.negri@unipv.it}}
%
%
\maketitle

\abstract{We consider the quasi-static evolution of a brittle layer on a stiff substrate; adhesion between layers is assumed to be elastic. 
Employing a phase-field approach we obtain the quasi-static evolution as the limit of time-discrete evolutions computed  by an alternate minimization scheme. 
We study the limit evolution, providing a qualitative discussion of its behaviour and 
a rigorous characterization, in terms of parametrized balanced viscosity evolutions.
Further, we study the transition layer of the phase-field, in a simplified setting, and show that it governs the spacing of cracks in the first stages of the evolution. Numerical results show a good consistency with the theoretical study and the local morphology of real life craquelure patterns.}

\section{Introduction}

\begin{figure}[h!]
\begin{center} \hspace{4pt}
\includegraphics[scale=0.2]{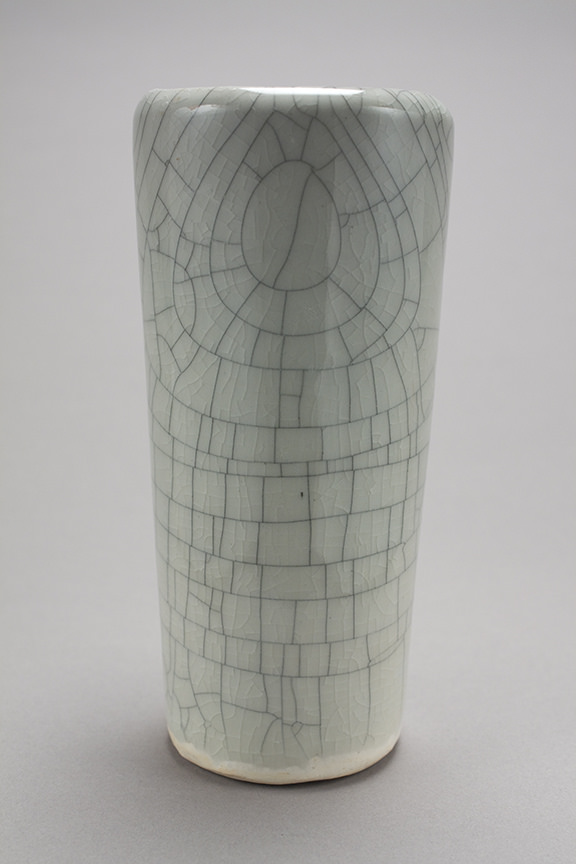}
 \quad 
		\includegraphics[scale=0.274]{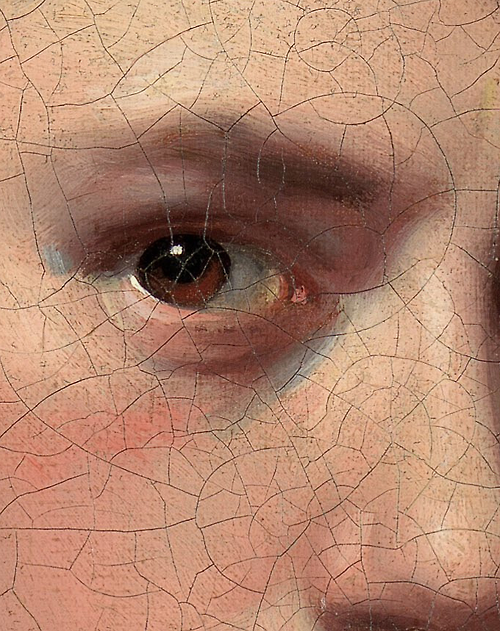} \break 
\caption{\label{f.iga} A regular pattern and a disordered network of cracks.}\end{center}
\end{figure}

Craquelure and crazing usually denote network of cracks which appear on thin superficial layers of materials (see Figure \ref{f.iga}). As a mathematical model for this kind of phenomena, we consider an elastic, brittle layer placed on a rigid adhesive substrate, which displaces the layer (and thus plays the role of a driving force). Having in mind the formation of patterns in a long time scale, we consider a quasi-static evolution driven by a time depending displacement of the substrate; in particular, we neglect diffusion of temperature and inertia. We employ a phase-field approach and we focus on a couple of complementary aspects: 
\begin{itemize}
\item characterization of quasi-static evolutions, obtained by time-discrete alternate minimization (staggered) schemes,
\item pattern formation and crack spacing as the result of different generations of cracks.
\end{itemize}
Other interesting aspects have been recently studied for similar mechanical systems, e.g.~\cite{BaldelliAndresBourdinMarigoMaurini_CMT13, AlessiFreddi_CS19, BraidesCausinSolci_PRSA}.

\begin{figure}[h!]
\begin{center}
	\includegraphics[scale=0.66]{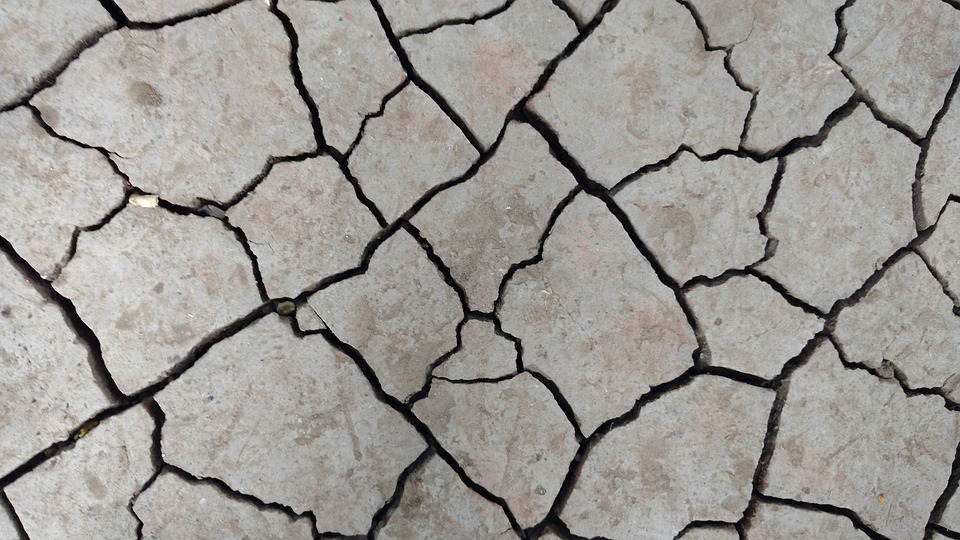} 
\quad %
		\includegraphics[scale=0.9]{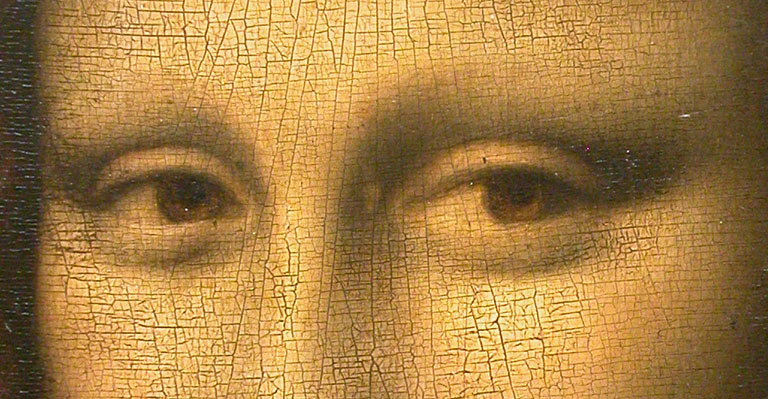} \qquad 
\end{center}
\caption{Other newtworks of cracks with different morphologies. }
\end{figure}

\separe

We assume that the reference configuration of the brittle film is an open bounded set $\Omega \subset \mathbb{R}^2$. The physical parameters are the Lam\`e coefficients $\lambda>0$ and $\mu>0$, the fracture toughness $G_c>0$ and the adhesion parameter $\beta>0$ (which measures the elastic response of the adhesive between the brittle layer and the substrate). 
Finally, we consider  an in-plane displacement $g$ which gives, as a function of time, the displacement of the substrate (note that here $g$ is a datum). For sake of simplicity we assume that $g(t)$ is of the form $t \hat{g} $ for some $\hat{g} \in H^1(\Omega, \R^2)$. 


The sharp crack energy associated to the system is, in some sense, a combination of Griffith energy \cite{Griff20}  and Mumford-Shah functional \cite{MumfordShah_CPAM89}; in the (weak) setting of $SBD^2$ spaces  
it reads
\begin{equation} \label{e.F}
\F ( t , u ) =  \tfrac12 \int_{\Omega \setminus J_u} W ( \strain(u)  ) \, dx + G_c \mathcal{H}^1 (J_u) + \beta \int_\Omega  | u - g (t) |^2 \, dx 
\end{equation}
where $W$ (depending on $\lambda$ and $\mu$) denotes the linear elastic energy density.
In this framework the set of discontinuity points $J_u$ represents the crack. 

Actually, in view of the numerical simulations, it is far more convenient to work with the phase-field  approximation of $\F$ given by \cite{BourdFrancMar00}
\begin{align}
\F_\eps ( t, u ,v) = 
\tfrac12 \int_\Omega & (v^2+ \eta_\eps ) W( \strain(u) ) \, dx \, + \nonumber \\ & + G_c \, \tfrac12 \int_\Omega  \eps^{-1} (v-1)^2   +   \eps | \nabla v |^2 \, dx + \beta \int_\Omega | u - g(t) |^2 \, dx .
\end{align}
It well known \cite{Chamb03} that the $\Gamma$-limit of the phase-field energy $\F_\eps (t, \cdot ,\cdot)$ is indeed the sharp crack energy $\F(t, \cdot)$. 

Now, let us introduce the time-discrete evolution which is used in the numerical simulations. Let $t_k = k \tau$ for $\tau>0$. Known the configuration $(u_{k-1}, v_{k-1})$ at time $t_{k-1}$ the updated configuration $(u_{k}, v_{k})$ at time $t_{k}$ is given by $v_k = \lim_{m \to +\infty} v_{k,m}$ and \ $u _k = \lim_{m \to +\infty} u_{k,m}$, where $u_{k,m}$ and $v_{k,m}$ are the steps of the staggered scheme:
$$\begin{cases}
	u_{k,m} \in \mbox{argmin} \, \big\{  \F_\epsilon  ( t_k \,, u  , v_{k,m-1} ) \big\} \phantom{\int} 
\\ 
          v_{k,m}  \in \mbox{argmin} \, \big\{  \F_\epsilon ( t_k \, , u_{k,m}  \, , v \, ) : v \le v_{k-1}  \big\} . \phantom{\int}
\end{cases} 
$$
The theoretical question we are interested is the characterization of the time continuous limit. 
A precise statement is contained in \S \ref{s.42}, asuuming a stronger irreversibility condition, i.e.~$v \le v_{k,m-1}$. Roughly speaking we will
 get an evolution $t \mapsto (u(t), v(t))$ which is possibly discontinuous and such that $(u(t), v(t))$ is a critical point of the energy, for every continuity point $t$. The technically difficult part, is instead the characterization of the behaviour in the discontinuity points, where the system makes an instantaneous (catastrophic) transition between $(u(t^-), v(t^-))$ and $(u(t^+), v(t^+))$ following some sort of staggered scheme. 

We finally provide some numerical results, which shows the pattern formation for a prototype problem and for the real life specimen of Figure \ref{f.iga}. 
It turns out that the pattern is formed by generations of cracks (which nucleate at different time) whose spacing follows a regular scheme: in the first generation, spacing of cracks is not predicted by minimizers (unless the domain is very short), it is instead dictated by the size of the boundary layer of the strain and of the phase-field function; further generations of cracks follow instead the dyadic behaviour of minimizers.


\section{Quasi-static evolutions by alternate minimization}

\subsection{Phase-field energy}

First of all let us introduce more rigorously the phase-field setting: for $\eps>0$ and $\eta_\eps>0$ we define the family of energies $\F_\eps : [0,T] \times H^1(\Omega; \R^2) \times H^1(\Omega, [0,1]) \to [0,+\infty)$ given by 
\begin{align}
\F_\eps ( t, u ,v) = 
\tfrac12 \int_\Omega & (v^2+ \eta_\eps ) W( \strain(u) ) \, dx \, + \nonumber \\ & + G_c \, \tfrac12 \int_\Omega  \eps^{-1} (v-1)^2   +   \eps | \nabla v |^2 \, dx + \beta \int_\Omega | u - g(t) |^2 \, dx .
\end{align}
Note that $\F_\eps ( t, \cdot, \cdot)$ is separately quadratic, i.e.~$\F_\eps ( t, u, \cdot)$ and $\F_\eps ( t, \cdot, v)$ are quadratic (positive) functionals; hence $\F_\eps ( t, u, \cdot)$ and $\F_\eps ( t, \cdot, v)$ are coercive and convex. On the other hand, $\F_\eps (t, \cdot, \cdot)$ is not (jointly) convex, however it is weakly lower semi-continuous in $H^1(\Omega; \R^2) \times H^1(\Omega, [0,1])$ (see e.g. \cite[Lemma 2.1]{Negri_ACV19}). In particular, for each time $t$ there exists a minimizer of the energy $\F_\eps ( t , \cdot,\cdot)$.

\separe

Now, let us turn to the $\Gamma$-limit of $\F_\eps ( t, \cdot, \cdot)$ as $\eps \to 0^+$. To this end, it is convenient to define the extended functional $\tilde\F_\eps : [0,T] \times L^2 (\Omega; \R^2) \times L^2 (\Omega, [0,1]) \to [0,+\infty]$ as
$$
	\tilde\F_\eps  (t, u ,v ) = 
	\begin{cases}
		\F _\eps (t,u,v) & \text{if $u \in H^1(\Omega; \R^2)$ and $v \in H^1(\Omega, [0,1])$,} \\
		+ \infty & \text{otherwise.} 
	\end{cases}
$$	
Then, by \cite{Chamb03} we known that the $\Gamma$-limit of $\F_\epsilon (t, \cdot, \cdot)$, with respect to the topology of $L^2 (\Omega; \R^2) \times L^2 (\Omega), $  is the functional 
$\tilde{F}(t, \cdot, \cdot) : [0,T] \times L^2 (\Omega; \R^2) \times L^2 (\Omega, [0,1]) \to [0,+\infty]$ given by 
$$
	\tilde\F  (t, u ,v ) = 
	\begin{cases}
		\F (t,u) & \text{if $u \in SBD^2 (\Omega)$ and $ v =1$ a.e.~in $\Omega$,} \\
		+ \infty & \text{otherwise.} 
	\end{cases}
$$
Let $(u_\eps, v_\eps)$ be a family of minimizers of $\F_\eps (t, \cdot,\cdot)$. Note that $u_\eps$ is bounded in $L^2$, then, by a fundamental result in the theory of $\Gamma$-convergence, see e.g.~\cite[Theorem 3.3]{Braides98}, we know that $u_\eps$ converge to a minimizer $u$ of $\F (t, \cdot)$; more precisely, there exists a subsequence converging to $u$ in the topology of $L^2 (\Omega; \R^2)$. On the contrary, not much is known on the convergence of the critical points of $\F_\eps(t, \cdot, \cdot)$ which will appear in the evolution. 

\subsection{Time-discrete evolution}

In this section we will describe the time-discrete evolution on which the numerical calculations are based. Here we can give only a brief description in the simplest possible setting, for complete proofs and generalizations the reader should make reference to \cite{KneesNegri_M3AS17, AlmiNegri_19}. 

For $\tau>0$ consider a time discretization $t_k = k \tau$. After setting the initial conditions $u_0$  and $v_0$, at time $t=0$, the configuration $(u_k, v_k)$, at time $t_k$, is computed from $(u_{k-1}, v_{k-1})$ by the following incremental scheme, known as alternate minimization algorithm. 
Let us introduce the further sequences $u_{k,m}$ and $v_{k,m}$, with $v_{k,0} = v _{k-1}$ and $u_{k,0} = u _{k-1}$. Then, for $m \ge 1$ we define
\begin{equation} \label{e.scheme1} \begin{cases}
	u_{k,m} \in \mbox{argmin} \, \big\{  \F_\epsilon  ( t_k \,, u  , v_{k,m-1} ) \big\} \phantom{\int} 
\\ 
          v_{k,m}  \in \mbox{argmin} \, \big\{  \F_\epsilon ( t_k \, , u_{k,m}  \, , v \, ) : v \le v_{k-1}  \big\} \phantom{\int}
\end{cases} 
\end{equation}
Note that this scheme takes full advantage of the fact that $\F ( t_k , \cdot , \cdot)$ is separately quadratic, and indeed its numerical implementation is very convenient, even if in practice the algorithm may converge quite slowly.
Then, (up to subsequences) we let $v_k = \lim_{m \to +\infty} v_{k,m}$ and \ $u _k = \lim_{m \to +\infty} u_{k,m}$. Note that the distance between $v_{k-1}$ and $v_k$ could be either small or large, in other terms, there is no a priori control on the speed $(v_{k} - v_{k-1})/\tau$ since the system is rate-independent. Moreover, note that $( u_k , v_k )$ is an equilibrium point for $\F_\epsilon (t_k , \cdot , \cdot)$, unilateral w.r.t. $v$. 

Finally, it is important to comment on the constraint $v \le v_{k-1}$ which models irreversibility. 
The convergence analysis under the latter constraint in the ''genuinely''  rate-independent setting is still open. At the current stage a couple of alternative are feasible. The first consists in solving a discrete viscous parabolic evolution for the phase-field, with the contraint $v \le v_{k-1}$, followed by a vanishing viscosity procedure: by \cite{AlmiBelzNegri_M2AN19} and \cite{Almi_19} this approach provides in the end a quasi-static evolution. The second option consists instead in replacing the constraint $v \le v_{k-1}$ by the (stronger) constraint $v \le v_{k,m-1}$. In this way monotonicity is imposed at each alternate iteration. On the theoretical level, this assumption allows to give a full characterization of the evolution as the time step $\tau \to 0$, without passing through viscosity solutions. On the other hand, in some cases this assumption may be too strong in the numerical simulations; we will see an example in \S \ref{s.irrrev}. Here we will follow the latter strategy.

\subsection{Time-continuous evolution \label{s.42}}

The time discrete scheme of the previous subsection gives for every $\tau>0$ a finite sequence $(v_k, u_k)$ in the points $t_k$. We will denote by $(u_\tau, v_\tau) : [0,T] \to H^1 (\Omega ; \R^2) \times H^1(\Omega)$ an interpolation of $(u_k, v_k)$ in the points $t_k$ (here we will enter into the delicate technical issue about the choice of the interpolation). Our goal is to characterize the limit of $(u_\tau, v_\tau)$ as $\tau \to 0$. 



\separe

First, we give a ``qualitative'' description of the limit evolution, as a function of time. We have already remarked that in this setting there is no a priori control of the length of $v_\tau(t_k) - v_\tau (t_{k-1})$ for $t_k = k \tau$. Indeed, fix $t \in [0,T]$ and let $k_\tau$ s.t. $t \in [k_\tau \tau , (k_\tau+1) \tau ]$. As $\tau \to 0$ it is clear that $k_\tau \tau \to t$ and $(k_\tau +1) \tau \to t$; on the contrary it may happen that $v_\tau (  k_\tau \tau ) \to v (t^-)$ and $v_\tau (  (k_\tau+1) \tau ) \to v(t^+)$ where $v(t^-) \neq v(t^+)$. In other terms the limit evolution $t \mapsto v(t)$ could be discontinuous in time. Note that this is observed in the numerical simulations. For this reason, in the limit as $\tau \to 0$ we expect an  evolution $ t \mapsto (u(t), v(t))$ of class $BV$ in time. Now, let us briefly describe the behaviour in continuity and discontinuity points (for more details we refer to \cite{AlmiNegri_19, KneesNegri_M3AS17}).

If $t \in [0,T]$ is a continuity point for the evolution then $(u(t), v(t))$ is an equilibrium configuration for the system, i.e. 
\begin{align*}
\partial_u \F_\epsilon ( t , u(t) , v(t) ) [\phi] & = 0  \ \text{ for every } \ \phi \in H^1(\Omega , \R^2) ,  \\
\partial_v \F_\epsilon ( t , u(t) , v(t) ) [\xi] & = 0  \ \text{ for every } \ \xi \in H^1(\Omega) \text{ with }  \xi \le 0 .
\end{align*}
Equilibrium of $v$ is unilateral because, by irreversibility, only negative variations are allowed. Note that the configuration $(u(t), v(t))$ is not necessarily a (joint) global or a local minimizer of $\F_\eps (t , \cdot, \cdot)$, it is actually a separate minimizer, by the separate convexity of $\F_\eps( t, \cdot, \cdot)$. 

On the contrary, if $t$ is a discontinuity point we expect an alternate ``evolution'', connecting $(u (t^-) , v(t^-))$ and $(u (t^+) , v(t^+))$. More precisely, in the limit the path between $(u (t^-) , v(t^-))$ and $(u (t^+) , v(t^+))$  is made of infinitely many intermediate configurations $ v_j \nearrow v (t^-)$ and $u_j \to u (t^-) $ connected by a (reverse) alternate scheme: i.e.
$$ \begin{cases}
	u _j \in \mbox{argmin} \, \big\{  \F_\epsilon  ( t \,, u  , v_{j+1} ) \big\} \phantom{\int} 
\\ 
          v_j  \in \mbox{argmin} \, \big\{  \F_\epsilon ( t \, , u_{j}  \, , v \, ) : v \le v_{j+1}  \big\} . \phantom{\int}
\end{cases} 
$$
In particular the instantaneous transition between $(u (t^-) , v(t^-))$ and $(u (t^+) , v(t^+))$ is not simultaneous in $u$ and $v$. Once again, this is confirmed by numerical results. 

\separe

In order to give a rigorous description of the limit evolution, it is necessary to introduce the derivatives of the energy $\F_\eps$, which will provide the driving forces for the evolution. 
We follow in particular \cite{AlmiNegri_19}. 
Clearly, we can take the partial derivatives of $\F_\eps$; moreover, we can define the following slopes
\begin{gather*}
|\partial_{u}\F_\epsilon|_{H^1}(t ,u ,v ) = \sup \{ - \partial_u \F_\eps  (t, u ,v) [ \phi ] : \phi \in H^1 (\Omega ; \R^2) \} ,
\\
|\partial_{z}^{-}\F_\epsilon|_{L^2}(t ,u ,v ) = \sup \{ - \partial_v \F_\eps  (t, u ,v) [ \xi ] : \xi \in H^1 (\Omega) \cap L^\infty(\Omega) \,, \xi \le 0 \} .
\end{gather*}
Note that the slope w.r.t.~$v$ is unilateral since, by irreversibility, we are interested only in negative variations. From the technical point of view it is fundamental that the slopes are weakly lower semicontinuous. Technically, we characterize the evolution in terms of a parametrized ``balanced viscosity solution'' \cite{MielkeRossiSavare_COCV12}. First, in order to describe both the behaviour in continuity and discontinuity points we employ an arc-length parametrization $s \mapsto ( t(s) , u(s) , v(s))$. In this setting, roughly speaking, discontinuity points (in time) correspond to intervals $[s^-, s^+]$ where $t$ is constant, and thus $t'= 0$; vice-versa, continuity points (in time) correspond to points $s$ where $t' (s) >0$. Then, the evolution is characterized by the following set of conditions: 

\begin{itemize}
\item [$(a)$] \ the map $t\colon [0,S]\to[0,T]$ is non-decreasing and surjective, 

\smallskip
\item [$(b)$] \ $(t,u,v)\in W^{1,\infty}([0,S];[0,T] \times H^1 ( \Omega ; \mathbb{R}^2) \times L^{2}(\Omega))$ with
\begin{displaymath}
|t'(s)|+\|u'(s)\|_{H^{1}}+\|v'(s)\|_{L^2}\leq 1\,,
\end{displaymath}

\item [$(c)$] \ $v$ is non-increasing and takes values in $[0,1]$,

\smallskip
\item [$(d)$] \ if $t'(s) >0$ then 
\begin{displaymath}
|\partial_{u}\F_\epsilon|_{H^1}(t(s),u(s),v(s))=0\qquad\text{and}\qquad|\partial_{z}^{-}\F_\epsilon|_{L^2}(t(s),u(s),v(s))=0 ,
\end{displaymath}

\item [$(e)$] \ for every $s\in[0,S]$ the following energy identity holds:
\begin{equation*}\label{e.eneq} 
\begin{split}
\F_\epsilon(t(s),& \ u(s),v(s)) = \F_\epsilon(0,u_0,v_0) +\int_{0}^{s} \!\! \partial_t \F_\eps ( t(r),u(r), v(r)) \,t'(r)\,dr \\
&-\int_{0}^{s} \!\! |\partial_{z}^{-}\F_\epsilon|_{L^2}(t(r),u(r),v(r))\, \| v' (r) \|_{L^2} \, dr \\ & 
-\!\int_{0}^{s}\! |\partial_{u}\F_\epsilon|_{H^1}(t(r),u(r),v(r))\, \| u' (s) \|_{H^1} \,dr  .
\end{split}
\end{equation*}
\end{itemize}
For further details the reader may refer to \cite{AlmiNegri_19}.

\section{A one-dimensional case study \label{s.1d}}

As a preliminary study, it is interesting to focus on the behaviour and on the pattern generated by the minimizers of the sharp-crack energy, even if the evolution follows only partially this scheme. To this end it is useful to consider a simpler one-dimensional setting, 
where $\Omega = (-L,L)$ and $g(t,x ) = t x $; the corresponding energy is then of the form
\begin{equation} \label{e.Fone}
   F (t , u ) =  \tfrac12 \int_{(-L,L)} \mu | u' |^2 \, dx + G_c \# (J_u) + \beta \int_{(-L,L)} | u - g(t) |^2 \, dx .
\end{equation}
We will assume that the initial configuration (at time $t=0$) is $g (0) = 0$.



\subsection{Continuous displacement and boundary layer \label{s.bl} }

\begin{figure}[h!]\begin{center}
\includegraphics[scale=0.65]{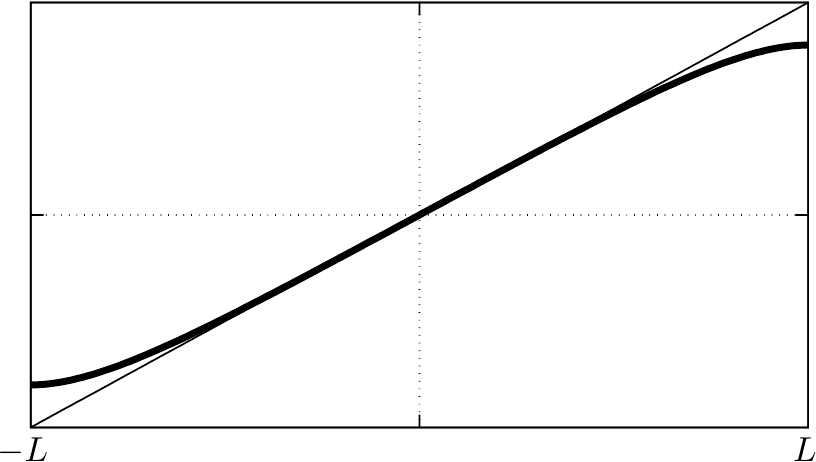} \qquad \includegraphics[scale=0.65]{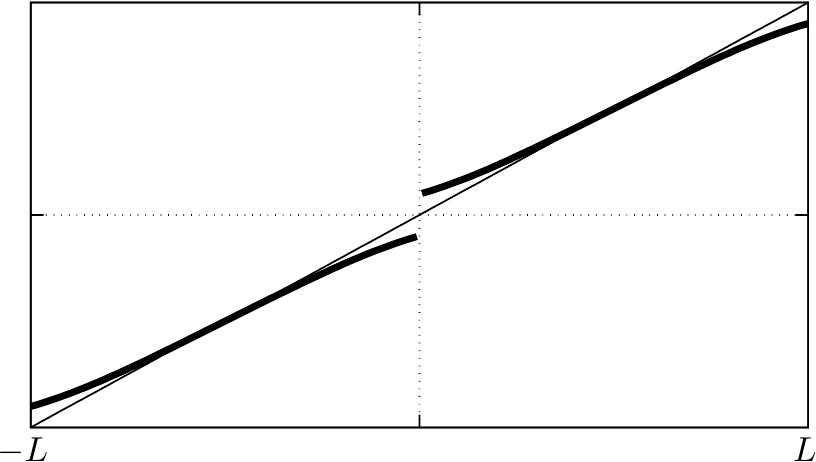}
\caption{\label{fig.u} Plot of $u_0 = u_L$ and $u_1$.}
\end{center}\end{figure}

First, let us consider a minimizer without cracks, i.e.~with $J_u = \emptyset$. In this case, the configuration  solves the Euler-Lagrange equation $- \mu u'' + 2 \beta u = 2 \beta g (t)$ in $(-L,L)$ with homogeneous Neumann boundary condition. 
The explicit solution is given by (see Figure \ref{fig.u})
$$	u_{0} (t,x) =  - t a_L \sinh ( \lambda x) + t x , \qquad \text{for $ \lambda = (2\beta / \mu)^{1/2}$ 
and $a_L = \frac{1}{\lambda \cosh(\lambda L)}$. }			$$
Then, a simple (but quite long) computation shows that 
$$
	F ( t , u_{0} (t)) = t^2 \mu \left[ L - \tfrac{1}{\lambda} \tanh(\lambda L)  \right] = t^2 \hat{F} (L) .
$$

\begin{figure}[h!]\begin{center}

\includegraphics[scale=0.64]{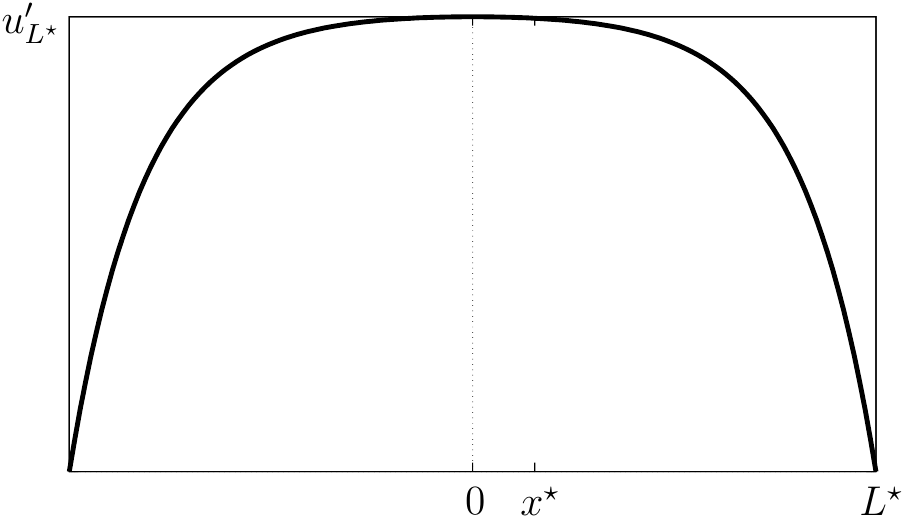}
\includegraphics[scale=0.64]{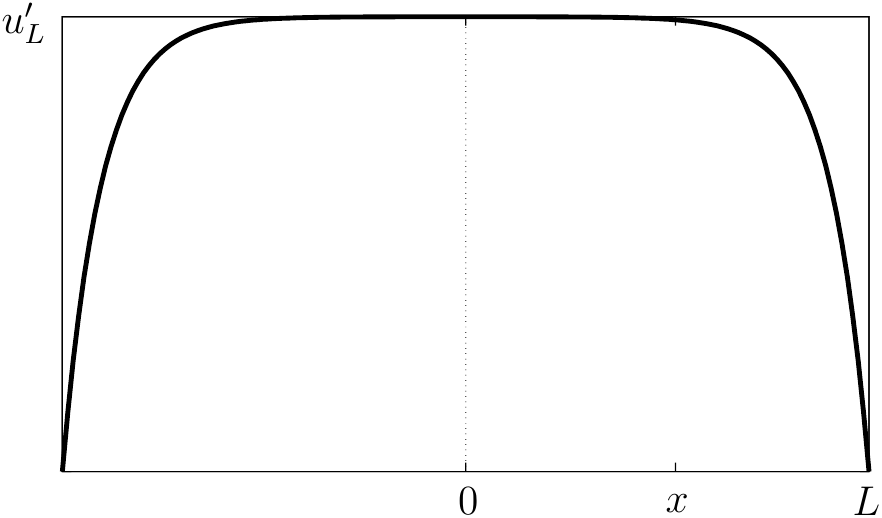}
\caption{\label{fig.bl} Plot of $u_{L^\star}$ for $L^\star=6.5$ and plot of $u_L$ for $L=12.5$.}
\end{center}\end{figure}

In the sequel it will be important we understand the dependence of $u_0$ on $L$. To this end, let us denote $u_L (x) = - a_L \sinh (\lambda x) + x$ the continuous solution $u_0$ in the interval $(-L,L)$ for $t=1$. Plotting the derivative 
$$	u'_L(x) = - a_L \lambda \cosh(\lambda x) + 1 = 1 - \frac{\cosh(\lambda x)}{\cosh(\lambda L)}  $$
for different values of $L$ (see Figure \ref{fig.bl}) shows that $u'_L$ changes in the boundary layers while it is almost constant in the ``interior''.

\begin{figure}[h!]\begin{center}
\includegraphics[scale=0.65]{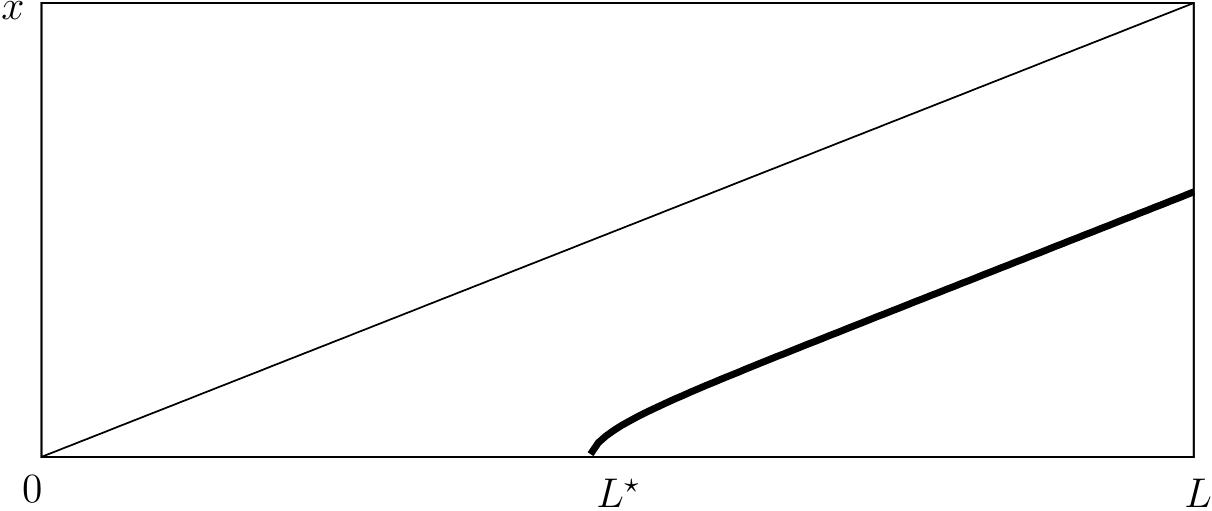}
\caption{\label{fig.lay} Plot of $L$ and $x$ (bold) from \eqref{e.xxstar} as a function of $L$. Values are  computed for $L^\star=6.5$ and $x^\star=0.5$. In this case we have $L - x \sim 5.5$.}
\end{center}\end{figure}

In order to provide some quantitative estimate of the size of the layer it is convenient to study how the derivative $u'_L(x)$
scales with $L$. Thus, given a point $x^\star \in (0,L^\star)$ let us compute the (corresponding) point $x \in (0,L)$ such that $u'_{L} (x) = u'_{L^\star} ( x^\star )$. The latter identity reads
\begin{equation} \label{e.xxstar}
	 \cosh( \lambda x ) = \frac{\cosh ( \lambda L) \cosh ( \lambda x^\star)}{\cosh (\lambda L^\star)} .
\end{equation}
A plot of the behaviour of $x$ as a function of $L$ (see Figure \ref{fig.lay}) shows that the distance between $x$ and $L$ is almost independent of $L$. Thus, we can expect that the length of the boundary layer remains almost constant with $L$.  As we will see in the numerical simulations of \S\,\ref{s.num}, the ``scale invariance'' of the boundary layer explains the periodic spacing of the first generations of cracks better than the dyadic structure (explained in the next subsections), which occurs only for small values of the size $L$ (more precisely, when $L$ is comparable with the size of layer).

In order to better explain the fact that $L - x$ is almost constant, we can consider the case when $x$ and $L$ are  large enough to approximate $\cosh$ with $\exp/2$, then the above identity can be approximated by 
$$
	\exp (\lambda x) \sim \exp (\lambda L) \exp( - \lambda b) \qquad 
	\text{for $b = - \ln \left( \frac{\cosh (\lambda x^\star)}{\cosh (\lambda L^\star)} \right) > 0$, } 
$$
which gives $x \sim L - b$. Roughly speaking, the transition layer is of constant size $b$. 

\subsection{First generation of cracks by minimality}

The content of this section follows the study performed in \cite{BaldelliAndresBourdinMarigoMaurini_CMT13}.
In the case of a single crack let us assume (for the moment) that $J_u = \{ 0 \}$, so that the interval $(-L,L)$ splits into the subintervals $(-L,0)$ and $(0,L)$; the minimizer, denoted by $u_{1}(t)$, can then be computed using the previous Euler-Lagrange equation in the subintervals $(-L,0)$ and $(0,L)$ (see again Figure \ref{fig.u}).  Then, the energy reads
$$
	F ( t , u_{1} (t) ) =  G_c + t^2 \, 2 \hat{F} (L/2) .
$$

Now we will compare the energy $F ( t , u_{0} (t))$ of the continuous solution with the energy $F ( t , u_{1} (t) )$ of the discontinuous solution. To this end,  it is interesting to consider in the $(t,L)$ plane the set of ``critical transition times'', satisfying
\begin{equation} \label{e.critical}
	t^2 ( \hat{F}(L) - 2 \hat{F}(L/2) )  = G_c .
\end{equation}
It is easy to see that  the function 
$$\Delta_2 (L) = \hat{F}(L) - 2 \hat{F}(L/2)=\tfrac{\mu}{\lambda} \left[ 2 \tanh(\lambda L /2) - \tanh(\lambda L ) \right] $$ vanishes for $L=0$ and is increasing  for $L>0$; hence it is positive and, given $L$,  there exists a time 
\begin{equation} \label{e.tcr}
	t_L = \left( \frac{G_c}{\hat{F}(L) - 2 \hat{F}(L/2)} \right)^{1/2} = \left( \frac{G_c}{\Delta_2 (L)} \right)^{1/2} 
\end{equation}
for which the two energies coincide. Further, $ t^2 ( \hat{F}(L) - 2 \hat{F}(L/2) ) < G_c $ for $t < t_L$ while $ t^2 ( \hat{F}(L) - 2 \hat{F}(L/2) ) > G_c $ for $t > t_L$; in other terms,
$$ F ( t , u_{0} (t) ) < F ( t , u_{1} (t) )  \ \text{ for $t < t_L$,} \qquad F ( t , u_{0}(t) ) > F ( t , u_{1} (t) )  \text{ for  $t > t_L$,}  $$  which justifies the name ''critical transition time'' at length $L$. 


\separe

Before proceeding, let us check that for a single crack the least energy is always assumed when $J_u = \{ 0 \}$. For $L_1 + L_2 =L$, we have to check the energy inequality 
$$
	G_c + \,t^2 \!\mu \left[ L_1 \!- \!\tfrac{1}{\lambda} \tanh(\lambda L_1) \right] + \,t^2 \!\mu \left[ L_2 \!- \! \tfrac{1}{\lambda} \tanh(\lambda L_2) \right] > G_c + \, t^2 \!\mu \left[ L \!- \! \tfrac{2}{\lambda} \tanh(\lambda L/2) \right] .
$$
After some algebraic manipulations, the previous inequality boils down to $$\tfrac12 ( \tanh(\lambda L_1) + \tanh(\lambda L_2)) < \tanh(\lambda (L_1 + L_2)/2), $$ which is true by the strict concavity of $\tanh$.

\separe

Now, let us consider the general case of $m$ cracks (for $m\ge1$). By the previous symmetry argument it is not restrictive to assume that $J_u = \{ - L + k L / (m+1) : \text{for $1 \le k \le m$} \}$. Arguing as above, the energy for $m$ cracks reads 
$$
	F ( t , u_m (t)) = m G_c + t^2 (m+1) \hat{F} (L/(m+1)) .
$$
We want to compare $F(t, u_0 (t))$ and $F(t, u_m (t))$. This time the set of critical points is defined by
$$
	t^2 \left[ \hat{F}(L) - (m+1) \hat{F}(L/(m+1)) \right] = m G_c ,
$$
which gives the ``critical transition time''
$$
	t_m^2 = \frac{m G_c}{\hat{F}(L) - (m+1) \hat{F}(L/(m+1))} .
$$
Let us show that the sequence $t_m$ is monotone increasing w.r.t.~$m$, independently of $L$. To this end it is convenient to study the function
$$
	s \mapsto \frac{(s-1)}{s \tanh(\lambda L /s) - \tanh(\lambda L) }
$$
which coincide with $t^2_m$ (up to a positive multiplicative constant) when $s=m+1$. A simple calculation shows that this function is increasing and thus $t_m$ is monotone increasing as well. Note that $t_1= t_L$. 

As a consequence, if the evolution is driven by energy minimization, in the first generation it appears a single crack in the center of the bar $(-L,L)$. Indeed, for $t \in [0, t_1)$ we have $F ( t , u_{0}(t) ) < F ( t , u_{1}(t) ) $ and thus $u (t)= u_{0} (t)$. For $t > t_1$ we have $F ( t_1 , u_{0}(t) ) > F ( t_1 , u_{1}(t) ) $, moreover, being $t_m > t_1$ for every $m > 1$, we have $F ( t, u_{m} (t)) > F ( t , u_{0} (t)) $, at least for $t \sim t_1$. Hence, we expect that $u (t) = u_{1} (t) $ for $t \in ( t_1, t_2)$ where $t_2$ denotes the time when the second generation of crack will occur. 

\begin{figure}\begin{center}
\includegraphics[scale=0.75]{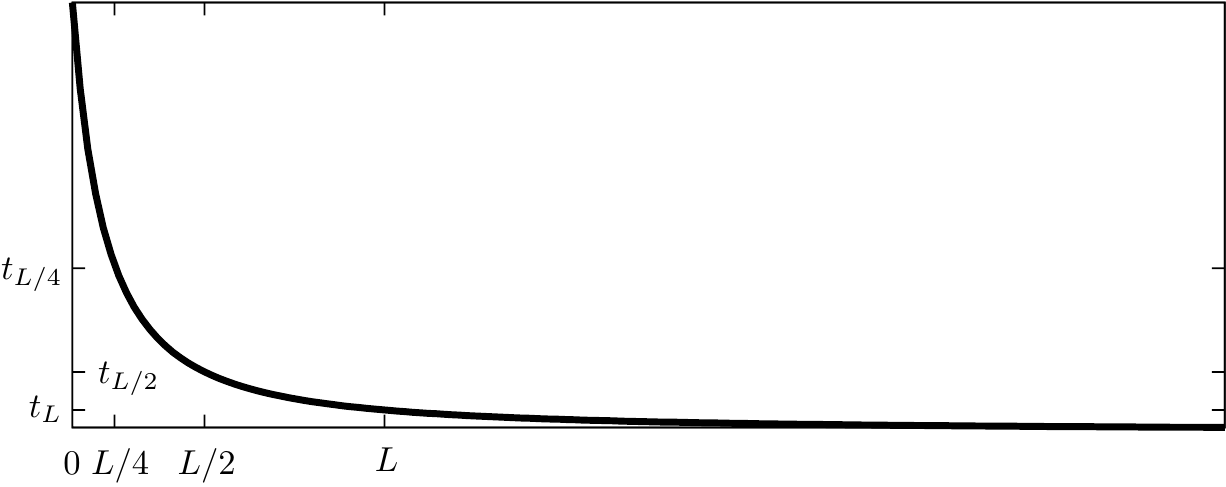}
\caption{\label{fig.tL} Plot of $t_L = (G_c / \Delta_2 (L))^{1/2}$.}
\end{center}\end{figure}

\subsection{Second and further generations of cracks by minimality}

When the first crack appears, at time $t_1$, the interval $(-L,L)$ is splitted into the subintervals $(-L,0)$ e $(0,L)$. By periodicity, this is equivalent to consider the behaviour of solutions in the interval $(-L/2,L/2)$.
Arguing exactly as in the previous section, just replacing $L$ with $L/2$, provides the existence of the critical transition time $t_{L/2}$ when the bar will split again. Note that $t_{L/2} > t_L$ because the function $t_L$ in \eqref{e.tcr} is decreasing w.r.t.~$L$ (remember that $\Delta_2$ is increasing with respect to $L$).

Proceeding by induction, we find the times $t_{L/2^k}$ (see Figure \ref{fig.tL}) when the dyadic structure evolves; in the physical literature this is often called ``halving'' effect.

\separe

To conclude this section, it is fair to remark that in this one dimensional example the transition between $u_{L}$ and $u_{L/2}$ cannot occur following a continuous energy decreasing path, because of the activation threshold $G_c$ which is payed as soon as the crack opens. In the numerical simulations, this topological problem is avoided by the phase-field regularization. Moreover, the uni-axial numerical results of \S \ref{s.num} we show that the first generation of cracks does not follow this dyadic structure (at least for $L$ large), rather the first crack pattern depends on the size of the boundary layer, described in the previous subsection.

\section{Alternate minimization in the one dimensional setting \label{s.uno-alt}}

In this section we briefly study the evolution generated by the alternate minimization scheme in the one-dimensional setting of \S \ref{s.1d}. For sake of simplicity we consider Dirichlet, instead of Neumann, boundary conditions. This example will be useful to understand the behaviour of the numerical results. We consider again the sharp crack energy $F$ of \eqref{e.Fone} and the corresponding phase-field energy 
\begin{align*}
F_\eps ( t, u ,v) = & 
\, \tfrac12 \int_{(-L,L)}  (v^2+ \eta_\eps ) \mu | u' |^2 \, dx \, + \nonumber \\ & + G_c \, \tfrac12 \int_{(-L,L)} \eps^{-1} (v-1)^2   +   \eps | v' |^2 \, dx + \beta \int_{(-L,L)} | u - g(t) |^2 \, dx .
\end{align*}
Here, we will further assume that $u(t, \pm L) = g ( t , \pm L) = \pm t L$. Fix $t_k = k \tau$ assume that $v_{k-1}=c_{k-1}$ is constant. We will show hereafter that the update $v_{k}$ is constant and that $u'_k$ is constant, as well. Looking at this result in the time discrete scheme, it turns out that for Dirichlet boundary conditions there will be no nucleation of cracks if the initial phase-field $v_0$ is homogeneous. Other interesting results on homogeneous states are contained in \cite{PhamMarigo_Je13}. 

For sake of simplicity, we will solve the following alternate scheme without the irreversibility constraint, which would actually not change the qualitative result,
\begin{equation*} \begin{cases}
	u_{k,m} \in \mbox{argmin} \, \big\{  F_\epsilon  ( t_k \,, u  , v_{k,m-1} ) : u (\pm L) =  g(\pm t_k L) \big\} \phantom{\int} 
\\ 
          v_{k,m}  \in \mbox{argmin} \, \big\{  F_\epsilon ( t_k \, , u_{k,m}  \, , v \, ) \big\} . \phantom{\int}
\end{cases} 
\end{equation*}
The Euler-Lagrange equation is of the form $- a_{k-1} u'' + 2 \beta u = 2 \beta g (t_k)$ in $(-L,L)$ with Dirichlet boundary conditions. 
The explicit solution is simply  $u_{k,m} (x) =  g (t_k,x) = - t_k x$ independently of $a_{k-1}$. Since $u_k$ is constant the Euler-Lagrange equation for $v$ is of the form $- v'' + b v = c $, whose solution $v_{k,m}$ is again a constant. In particular the staggered scheme finds a critical point after two iterations.


\section{Numerical results for uni-axial problems \label{s.num}}

In the following subsections are reported the numerical results obtained for a bar $\Omega = (-L , L) \times (-H,H)$ of variable length and fixed width $2H=5$. We have chosen fracture toughness $G_c=1.0$, Young modulus  $E=1.0$, Poisson ratio $\nu=0.15$, adhesive constant $\beta=0.15$. Moreover, we assume that the displacement of the substrate is of the form $g (t, x) = (t x_1, 0)$. Since $g$ is uni-axial we expect solutions in accordance with the theoretical arguments of \S\,\ref{s.1d}.   
Since the datum $g$ is monotone we neglected the irreversibility constraint in the alternate minimization scheme \eqref{e.scheme1}. With this choice crack patterns look sharper. A comparison between the solutions obtained with and without irreversibility constraint (see \S\,\ref{s.irrrev}) show that the crack patterns behave in a similar way. Numerical results have been computed using FreeFem++ \cite{Hecht_JNM12}. 

\subsection{Dyadic structure for short bars}

\begin{figure}[h!]\begin{center}
\includegraphics[scale=0.4]{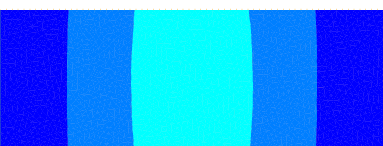} \qquad 
\includegraphics[scale=0.4]{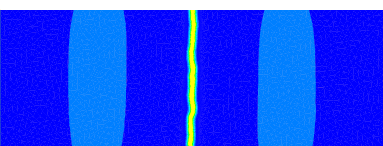} \\ 
\vspace{12pt}
\includegraphics[scale=0.4]{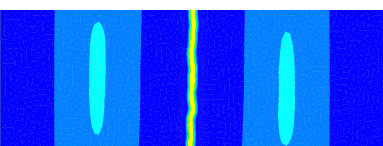} \qquad 
\includegraphics[scale=0.4]{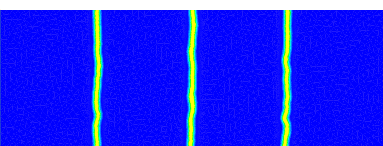}
\caption{\label{fig.rev} Phase field $v$ at time $t=\{3.0, 3.1,3.7,3.8\}$ computed without irreversibility constraint; blue corresponds to sound material, i.e.~$v=1$, while yellow corresponds to cracks, i.e.~$v=0$.}
\end{center}\end{figure}

In this subsection we present and discuss the numerical results obtained with $L=6.5$. The configuration of the phase field $v$ shown in Figure \ref{fig.rev} are computed with the unconstrained version of scheme \eqref{e.scheme1}, i.e. without any constraint on $v$.

First, it is interesting to comment on the very first image. 
In this case there are no pre-existing fractures, hence nucleation is a fundamental ingredient in the evolution. Figure \ref{fig.rev} shows that nucleation requires first a diffuse damage, with the phase-field $v$ decreasing in a wide region, and then concentration, in a second stage (see also \cite[Fig 4.8]{Kuhn_2013}). Note that, in the first image, the profile of $v$ is closely related to the profile of $u'_L$ (the derivative of the one dimensional solution  $u_L$) presented in \S  \ref{s.bl}. Indeed, since the datum is uniaxial and Poisson ratio $\nu$  is small, we can expect that $u'_L$ provides a good approximation of the strain energy density; then, in phase-field models the larger is the elastic energy density the smallest is the phase-field. This is why the profile of $v$ is consistent with that $u'_L$ (compare the first plot of Figure \ref{fig.bl} with the first image of Figure \ref{fig.rev}).

Then, the central crack nucleates and the domain splits. At that point, a similar evolution ``restarts'' in each subinterval, leading to the second generation of cracks.

Note that in this case the ``transition layer'' for $u'_L$ is, roughly speaking, of the same size of the half interval $(0,L)$. For this reason the evolution generates a dyadic fracture pattern, clearly visible in the last image of Figure \ref{fig.rev}. As we will see, the pattern of the first generation is not dyadic for larger $L$. 


\pagebreak

\subsection{Periodic patterns for longer bars \label{s.L12}}

In this subsection we present the numerical results obtained for $L=12.5$. In the second picture, when  the first generation of crack appears (at time $t=2.8$) it is evident that the position of the first crack  (from left to right)  corresponds to the transition region of $v$, which, comparing with Figure \ref{fig.bl}, corresponds to the transition layer of $u'_L$. Moreover, all the cracks, apart from the innermost, are equally spaced and their distance is approximately $5.5$, which is the value expected from the estimates of the size of the transition layer (see Figure \ref{fig.lay}). Note that, following the intermediate steps of the alternate minimization algorithm (at time $t=2.8$), the outermost cracks nucleates firts, followed by the innermost. 
The reason behind the fact that cracks do not nucleate in the center of the bar is not fully clear, however, it should be found in the fact that the phase field $v$ is almost constant in the inner part of the domain and thus its evolution does not promote concentration of strain and nucleation of cracks, as described in \S \ref{s.uno-alt}. The second generation of cracks follows instead the dyadic structure: this is due to the fact that the spacing of  the cracks after the first generation is comparable with the size of the transition layer, as in the previous example.

\begin{figure}[h!]\begin{center}
\includegraphics[scale=0.75]{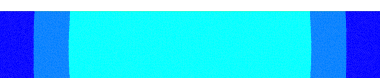} \\
\includegraphics[scale=0.75]{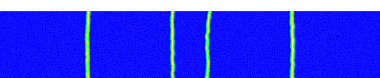} \\
\includegraphics[scale=0.75]{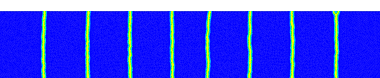} \\
\caption{\label{fig.long} Evolution for $L=12.5$ at time $t=2.8$, $t=3.1$ and $t=4.3$}
\end{center}\end{figure}

In conclusion, we could say that the first generation of cracks is determined by a characteristic length (the width of the transition layer) with cracks nucleating from the boundary to the interior, while the further generations evolve according to the dyadic scheme. Numerical results for longer bars confirm this behaviour. Finally we remark this type of evolution, triggered by the boundary layer, has been observed also in  numerical simulations of brittle layered materials, see \cite[Fig.~7]{AlessiFreddi_CS19}, and soil drying (which share several features with our problem), see \cite[Fig.~6]{VoPouyaHemmatiTang_CG17}.

\subsection{Irreversibility \label{s.irrrev}}

In this section we briefly report the numerical results obtained with the scheme \eqref{e.scheme1} under the constraint $v \le v_{k-1}$ in the case $L=6.5$. Clearly, the main difference is the fact that here the phase-field $v$ cannot increase after the nucleation of cracks, because of the irreversibility constraint. From the results, it is evident that this evolution is very similar to the evolution of Figure \ref{fig.rev}, at least as far as the structure of the crack pattern, the time and the nucleation sites.  

In the case $L=12.5$ the behaviour is similar to that of Figure \ref{fig.long}.


\begin{figure}[h!]\begin{center}
\includegraphics[scale=0.4]{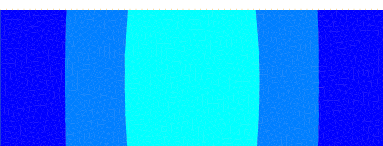} \qquad 
\includegraphics[scale=0.4]{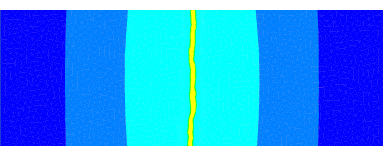} \\ 
\vspace{12pt}
\includegraphics[scale=0.4]{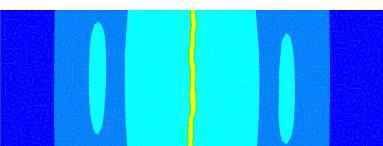} \qquad 
\includegraphics[scale=0.4]{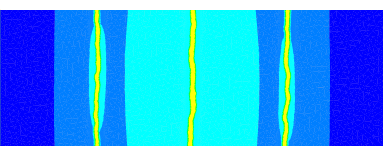}
\caption{\label{fig.irr} Phase field $v$ at time $t=\{3.0, 3.1,3.7,3.8\}$ computed with the time irreversibility constraint.}
\end{center}\end{figure}

%

\section{Local craquelure patterns for a real life specimen \label{l.real}}

\begin{figure}[ht!]

%
%
\begin{center}
\includegraphics[scale=0.6]{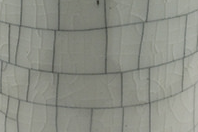} \quad \includegraphics[scale=0.256]{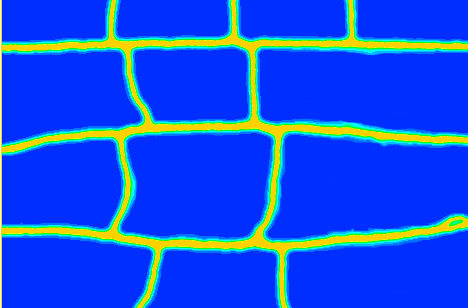} \qquad 
\includegraphics[scale=0.75]{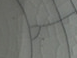} 

\vspace{-80pt} \hspace{244pt} \includegraphics[scale=0.25]{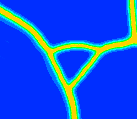}

\vspace{68pt}

\includegraphics[scale=0.5]{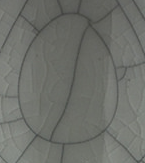} \quad \includegraphics[scale=0.34]{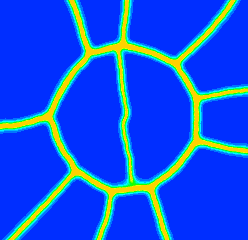} \qquad  
\includegraphics[scale=0.6]{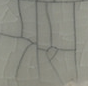} \quad 
%
\includegraphics[scale=0.19]{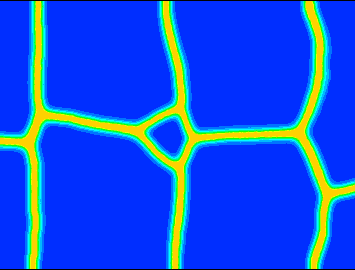} 


\caption{\label{fig-vas} Details of some numerical results obtained with different values of $\beta$ and $A$ compared with similar real life craquelure morphologies.}
\end{center}
\end{figure}

In order to produce some more realistic craquelure patterns and crack morphologies, we ran several numerical experiment changing both the value of the adhesive parameter $\beta$ and the datum $g$, of the form $g(x) = Ax$ for $A \in \R^{2 \times 2}$. 
Even if the values are not provided by experimental measurements, it is clear from Figure \ref{fig-vas} that 
this quasi-static phase-field model captures quite well the local morphological features of the patterns in Figure \ref{f.iga}. 
 Note that these plots do not show the entire evolution leading to the formation of the patterns but they show the final snapshot, at a certain time $T$. In the evolution the cracks actually appear at different times, following a scheme reminiscent of the uni-axial case.

In the bi-axial setting it is more difficult to explain the behaviour of the evolution. However, some observations are due as far as crack junctions, which of course did not occur in the uniaxial setting.
First, remember that for local minimizers of the (scalar) Mumford-Shah functional only Y-junctions (i.e. triple junctions) at $\tfrac23 \pi$ angles occur, see \cite[Theorem 2.1]{MumfordShah_CPAM89}. We may expect a similar behaviour in our linear elasticity context. However, in the first image only T-junctions occur, this is simply due to the fact that horizontal cracks nucleate in the first generation and vertical cracks in the second generation, therefore, once horizontal straight cracks are formed it is no longer possible to have $\tfrac23 \pi$ angles. In this example, it would be natural to see orthogonal crossing of cracks; on the contrary cracks seem to ``shift'' passing from one horizontal stripe to the other. This feaure could be due to remeshing and to the fact that horizontal cracks are not equally spaced, however, it is interesting to note that this feature occurs also in several points of the real life picture. In some cases  (see the pictures on the right) the crack behaves instead as a local minimizer, splitting a T- or an X-junction into three (or more) Y-junctions, with angles close to the optimal value $\tfrac23 \pi$, compare with \cite[Figure 14]{MumfordShah_CPAM89}. In this way, small triangular regions are formed. Once again, it is worth to remark that this the case both in the numerical experiment and in the real life specimen.



\bibliographystyle{plain}
\bibliography{references.bib} 

\end{document}